\documentclass[11pt]{amsart}
\usepackage[latin9]{inputenc}
\usepackage{amsthm}
\usepackage{amstext}
\usepackage{amssymb}

\makeatletter

\providecommand{\tabularnewline}{\\}

\numberwithin{equation}{section}
\numberwithin{figure}{section}
\theoremstyle{plain}
\newtheorem{thm}{\protect\theoremname}
  \theoremstyle{remark}
  \newtheorem{rem}[thm]{\protect\remarkname}

\title{}
\def\qed{\quad\vrule height4.17pt width4.17pt depth0pt}
\usepackage{amsfonts}
\usepackage{enumerate}

\newtheorem{theorem}{Theorem}
\theoremstyle{plain}

\newtheorem{definition}{Definition}
\newtheorem{example}{Example}

\newtheorem{lemma}{Lemma}

\newtheorem{remark}{Remark}

\numberwithin{equation}{section}

\makeatother

  \providecommand{\remarkname}{Remark}
\providecommand{\theoremname}{Theorem}

\begin{document}

\title{The Multiplication Group of an AG-group}

\author{M. Shah}

\email{shahmaths\_problem@hotmail.com}

\author{A. Ali}

\email{dr\_asif\_ali@hotmail.com}

\address{Department of Mathematics, Quaid-i-Azam University, Islamabad, Pakistan. }

\author{I. ahmad{*}}

\address{Department of Mathematics, University of Malakand, Chakdara Dir(L),
Pakistan.}

\email{iahmaad@hotmail.com}

\author{V. Sorge}

\email{V.Sorge@cs.bham.ac.uk}

\address{School of Computer Science, University of Birmingham, UK.}

\keywords{multiplication group, inner mapping group, AG-group, translations.\\
{*} Corresponding Author}

\begin{abstract}
We investigate the multiplication group of a special class of quasigroup
called AG-group. We prove some interesting results such as: the multiplication
group of an AG-group of order $n$ is non-abelian group of order $2n$
and its left section is an abelian group of order $n$. The inner
mapping group of an AG-group of any order is a cyclic group of order
$2$. 
\end{abstract}
\maketitle

\section{Introduction and Preliminaries}

A groupoid $G$ is an\emph{ AG-group} if it satisfies: $(i)$ $(xy)z=(zy)x,$
$\forall x,y,z\in G$. $(ii)$ There exists left identity $e\in G$
(that is, $ex=x,\forall x\in G$). $(iii)$ For every $x\in G$ there
exists $x^{^{-1}}\in G$ such that $x^{^{-1}}x=xx^{^{-1}}=e$. $x$
and $x^{^{-1}}$ are called inverses of each other. 

AG-group is a subclass of cancellative AG-groupoids \cite{SSA}. Some
basic properties of AG-groups have been derived in \cite{SA}, and
fuzzification of AG-groups can be seen in \cite{fuzz,fuzHt}. AG-group
is a generalization of abelian group and is a special quasigroup.
AG-groups have been counted computationally in \cite{SGS} and algebraically
in \cite{shahthesis}. The counting of AG-groups up to order $6$
can also be found in \cite{DSV}. AG-groups have been studied as a
generalization of abelian group as well as a special case of quasigroups
in \cite{shahthesis}. The present paper discovers the multiplication
group and inner mapping group of an AG-group. Multiplication group
and inner mapping group of a loop have been investigated in a number
of papers for example \cite{fixtwo,centnil,cycinnermap,contransfin,contrans,prodtwoprimes,dihedtwogroups,innermapabelian,transnonabelian}.
This has always been remained the most interesting topic of group
theorists in loop theory. Quasigroup does not have inner mapping group
because it does not have an identity element unless it is not a loop.
But an AG-group though not a loop but it has a left identity so it
has multiplication group as well as inner mapping group. We will prove
here some interesting results about the multiplication group and inner
mapping group of an AG-group that do not hold in case of a loop. For
example for an AG-group $G$ of order $n$ the left section, $L_{S}$
is an abelian group of order $n$. Its multiplication group is a nonabelian
group of order $2n$. The inner mapping group of an AG-group is always
a cyclic group of order $2$ regardless of its order. The following
lemma of \cite{SA} will be used in proofs .

\begin{lemma} \label{Pre-L1} Let $G$ be an AG-group. Let $a,b,c,d\in G$
and $e$ be the left identity in $G$. Then the following conditions
hold in $G$. 

\begin{enumerate}[(i)]

\item $(ab)(cd)=(ac)(bd)$ medial law; 

\item $ab=cd\Rightarrow ba=dc$;

\item $a\cdot bc=b\cdot ac$; 

\item $(ab)(cd)=(db)(ca)$ paramedial law; 

\item $(ab)(cd)=(dc)(ba)$; 

\item $ab=cd\Rightarrow d^{-1}b=ca^{-1}$; 

\item If $e$ the right identity in $G$ then it becomes left identity
in $G$, i.e, $ae=a\Rightarrow ea=a$;

\item $ab=e\Rightarrow ba=e$; 

\item $(ab)^{-1}=a^{-1}b^{-1}$; 

\item $a(b\cdot cd)=a(c\cdot bd)=b(a\cdot cd)=b(c\cdot ad)=c(a\cdot bd)=c(b\cdot ad)$; 

\item $a(bc\cdot d)=c(ba\cdot d)$; 

\item $(a\cdot bc)d=(a\cdot dc)b$; 

\item $(ab\cdot c)d=a(bc\cdot d)$. 

\end{enumerate}

\end{lemma}

\section{Multiplication Group of an AG-group}

Let $G$ be an AG-group and $a$ be an arbitrary element of $G$.
The mapping $L_{a}:G\rightarrow G$ defined by $L_{a}(x)=ax$ is called
left translation on $G$. Similarly the mapping $R_{a}:G\rightarrow G$
defined by $R_{a}(x)=xa$ is called right translation on $G$. 

Our first result establish the relation between a left translation
and a right translation.

\begin{lemma} \label{Ch7.Sec2.L2} Let $G$ be an AG-group. Let $a,b\in G$
and $e$ be the left identity in $G$. Then 

\begin{enumerate}[(i)]

\item $L_{a}R_{b}=R_{ab}.$ 

\item $R_{a}R_{b}=L_{ab}.$ 

\item $L_{a}L_{b}=R_{(ae)}R_{b}.$ 

\item $L_{a}L_{b}=L_{(ae)b}=L_{(be)a}.$ 

\item $R_{a}L_{b}=R_{(ae)b}.$ 

\item $L_{a}L_{b}=L_{b}L_{a}.$ 

\item $R_{a}L_{b}=R_{b}L_{a}.$ 

\end{enumerate}

\end{lemma}

\emph{Proof. }Let $a,b\in G$ and $e$ be the left identity in $G$.
Then 

\begin{enumerate}[(i)]

\item $L_{a}R_{b}(x)=L_{a}(xb)=a(xb)=x(ab)=R_{ab}(x)$ $\Rightarrow L_{a}R_{b}=R_{ab}$. 

\item $R_{a}R_{b}(x)=R_{a}(xb)=(xb)a=(ab)x=L_{ab}(x)$ $\Rightarrow R_{a}R_{b}=L_{ab}.$ 

\item $L_{a}L_{b}(x)=L_{a}(bx)=a(bx)=(ea)(bx)=(xb)(ae)=R_{(ae)}(xb)$
$=R_{(ae)}R_{b}(x)\Rightarrow L_{a}L_{b}=R_{(ae)}R_{b}.$ 

\item By (ii) and (iii) and left invertive law. 

\item $R_{a}L_{b}(x)=R_{a}(bx)=(bx)a=(bx)(ea)=(ae)(xb)=L_{ae}(xb)$
$=L_{ae}R_{b}(x)\Rightarrow R_{a}L_{b}=L_{ae}R_{b}\Rightarrow R_{a}L_{b}=R_{(ae)b}$,
by (i). 

\item $L_{a}L_{b}=L_{(be)a}$ $\Rightarrow L_{a}L_{b}=L_{b}L_{a}$,
by (iv). 

\item $R_{a}L_{b}=R_{(ae)b}$ $=R_{(be)a}=R_{b}L_{a}$, by left invertive
law and (v). $\qed$

\end{enumerate}
\begin{rem}
From Lemma \ref{Ch7.Sec2.L2} we note that if $G$ is an AG-group,
then the left translation $L_{a}$ and the right translation $R_{a}$
behave like an even permutation and an odd permutation respectively,
that is; 
\[
L_{a}L_{a}=L_{a},R_{a}R_{a}=L_{a},L_{a}R_{a}=R_{a},R_{a}L_{a}=R_{a}.
\]
 
\end{rem}
Next we recall the following definition. \begin{definition} Let $G$
be an AG-group. Then the set $L_{S}=\{L_{a}:L_{a}(x)=ax\,\forall$
$x\in G\}$ is called \textbf{left section} of $G$ and the set $R_{S}=\{R_{a}:R_{a}(x)=xa\,\forall$
$x\in G\}$ is called \textbf{right section} of $G.$ \end{definition}

\begin{definition} Let $G$ be an AG-group. Then the set $\langle L_{a},R_{a}:a\in G\rangle$
forms a group which is called multiplication group of the AG-group
$G$ and is denoted by $M(G)$ i.e $M(G)=\langle L_{a},R_{a}:a\in G\rangle.$
\end{definition}

\bigskip{}

We remark that left section of a loop is not a group but left section
of an AG-group does form a group as we prove it in the following theorem.

\begin{theorem} Let $G$ be an AG-group of order $n$. Then $L_{S}$
is an abelian group of order $n$. \end{theorem}

\begin{proof} By definition $L_{S}=\{L_{a}:L_{a}(x)=ax\,\forall$
$x\in G,a\in G\}.$ Let $L_{a},L_{b}\in L_{S}$ for some $a,b\in G.$
Then by Lemma \ref{Ch7.Sec2.L2} (iv), we have $L_{a}L_{b}=L_{(ae)b}\in L_{S}\Rightarrow L_{S}$
is an AG-groupiod. $L_{e}L_{a}=L_{(ee)a}=L_{a}$ and $L_{a}L_{e}=L_{(ae)e}=L_{(ee)a}=L_{a}.$
Therefore $L_{e}$ is the identity in $L_{S}.$

\noindent Let $L_{a},L_{b},L_{c}\in L_{S}.$ Then $(L_{a}L_{b})L_{c}=L_{(ae)b}L_{c}=L_{[\{(ae)b\}e]c}=L_{(ce)((ae)b)}=L_{(ce)((be)a)}=L_{(ae)((be)c)}$$=L_{a}L_{(be)c}=L_{a}(L_{b}L_{c}).$

\noindent Let $L_{a}\in L_{S}\Rightarrow a\in G\Rightarrow a^{-1}\in G\Rightarrow a^{-1}e\in G.$
Let $a^{-1}e=b$ then $L_{b}\in L_{S}.$ Now $L_{a}L_{b}=L_{(ae)b}=L_{(ae)(a^{-1}e)}=L_{e}=L_{b}L_{a}\Rightarrow L_{b}$
is the inverse of $L_{a}.$ Thus $L_{S}$ is a group. Since from Lemma
\ref{Ch7.Sec2.L2}, we have $L_{a}L_{b}=L_{b}L_{a}.$ Therefore $L_{S}$
is an abelian group. \end{proof}

\noindent We illustrate the above result by an example. \begin{example}
\label{Ch7.Sec2.E1} An AG-group of order $3:$

\begin{center}
$\begin{tabular}{l|lll}
 \ensuremath{\cdot} &  \ensuremath{0} &  \ensuremath{1} &  \ensuremath{2}\\
\hline \ensuremath{0} &  \ensuremath{0} &  \ensuremath{1} &  \ensuremath{2}\\
\ensuremath{1} &  \ensuremath{2} &  \ensuremath{0} &  \ensuremath{1}\\
\ensuremath{2} &  \ensuremath{1} &  \ensuremath{2} &  \ensuremath{0}
\end{tabular}$ 
\par\end{center}

\end{example}

The Multiplication group of the AG-group given in Example \ref{Ch7.Sec2.E1}
is isomorphic to $S_{3},$ the symmetric group of degree $3$ as the
following example shows.

\begin{example} \label{Ch7.Sec2.newE1} Multiplication group of the
AG-group given in Example \ref{Ch7.Sec2.E1}. 

\begin{center}
\begin{tabular}{l|llllll}
$\cdot$  & $L_{0}$  & $L_{1}$  & $L_{2}$  & $R_{0}$  & $R_{1}$  & $R_{2}$ \tabularnewline
\hline 
$L_{0}$  & $L_{0}$  & $L_{1}$  & $L_{2}$  & $R_{0}$  & $R_{1}$  & $R_{2}$ \tabularnewline
$L_{1}$  & $L_{1}$  & $L_{2}$  & $L_{0}$  & $R_{2}$  & $R_{0}$  & $R_{1}$ \tabularnewline
$L_{2}$  & $L_{2}$  & $L_{0}$  & $L_{1}$  & $R_{1}$  & $R_{2}$  & $R_{0}$ \tabularnewline
$R_{0}$  & $R_{0}$  & $R_{1}$  & $R_{2}$  & $L_{0}$  & $L_{1}$  & $L_{2}$ \tabularnewline
$R_{1}$  & $R_{1}$  & $R_{2}$  & $R_{0}$  & $L_{2}$  & $L_{0}$  & $L_{1}$ \tabularnewline
$R_{2}$  & $R_{2}$  & $R_{0}$  & $R_{1}$  & $L_{1}$  & $L_{2}$  & $L_{0}$ \tabularnewline
\end{tabular}
\par\end{center}

\end{example}

\noindent Here $L_{S}=\{L_{0},L_{1},L_{2}\}$ which is an abelian
group as the following table shows:

\begin{center}
$\bigskip\begin{tabular}{l|lll}
 \ensuremath{\cdot} &  \ensuremath{L_{0}} &  \ensuremath{L_{1}} &  \ensuremath{L_{2}}\\
\hline \ensuremath{L_{0}} &  \ensuremath{L_{0}} &  \ensuremath{L_{1}} &  \ensuremath{L_{2}}\\
\ensuremath{L_{1}} &  \ensuremath{L_{1}} &  \ensuremath{L_{2}} &  \ensuremath{L_{0}}\\
\ensuremath{L_{2}} &  \ensuremath{L_{2}} &  \ensuremath{L_{0}} &  \ensuremath{L_{1}}
\end{tabular}$ 
\par\end{center}

\noindent But $R_{S}=\{R_{0},R_{1},R_{2}\}$ does not form an AG-group
as the following table shows:

\begin{center}
\begin{tabular}{l|lll}
$\cdot$  & $R_{0}$  & $R_{1}$  & $R_{2}$ \tabularnewline
\hline 
$R_{0}$  & $L_{0}$  & $L_{1}$  & $L_{2}$ \tabularnewline
$R_{1}$  & $L_{2}$  & $L_{0}$  & $L_{1}$ \tabularnewline
$R_{2}$  & $L_{1}$  & $L_{2}$  & $L_{0}$ \tabularnewline
\end{tabular}
\par\end{center}

\begin{remark} Right section does not form even an AG-groupoid. \end{remark}

Lemma \ref{Ch7.Sec2.L2} guarantees that for an AG-group $G,$ $M(G)=\langle L_{a},R_{a}:a\in G\rangle=\{L_{a},R_{a}:a\in G\}$.

\begin{theorem} Let $G$ be an AG-group of order $n$. The set $\{L_{a},R_{a}:a\in G\}$
forms a non-abelian group of order $2n$ which is called multiplication
group of the AG-group $G$ and is denoted by $M(G)$ i.e $M(G)=\{L_{a},R_{a}:a\in G\}$.
\end{theorem}

\begin{proof} From Lemma \ref{Ch7.Sec2.L2}, it is clear that $M(G)$
is closed. $L_{e}$ plays the role of identity as $L_{a}L_{e}=L_{e}L_{a}=L_{a}$.
\[
R_{a}L_{e}=R_{(ae)e}=R_{(ee)a}=R_{a}=R_{ea}=L_{e}R_{a}.
\]

\noindent Let $L_{a}\in M(G)\Rightarrow a\in G\Rightarrow a^{-1}\in G\Rightarrow R_{a^{-1}}\in M(G)$
and $R_{a}R_{a^{-1}}=L_{aa^{-1}}=L_{e}=L_{a^{-1}a}=R_{a^{-1}}R_{a}.$
Therefore $R_{a^{-1}}$ is the inverse of $R_{a}$ in $M(G).$ Associativity
in $M(G)$ follows from the associativity of mappings$.$ Thus $M(G)$
is a group. Note that $M(G)$ is non-abelian because $R_{a}R_{b}\neq R_{b}R_{a}$
by \ref{Ch7.Sec2.L2} (ii). \end{proof}

\bigskip{}

To make things a bit more clearer we consider the following examples.

\begin{example} \label{Ch7.Sec2.E2} An AG-group of order $4$. 

\begin{center}
\begin{tabular}{l|llll}
$\cdot$  & $0$  & $1$  & $2$  & $3$ \tabularnewline
\hline 
$0$  & $0$  & $1$  & $2$  & $3$ \tabularnewline
$1$  & $1$  & $0$  & $3$  & $2$ \tabularnewline
$2$  & $3$  & $2$  & $1$  & $0$ \tabularnewline
$3$  & $2$  & $3$  & $0$  & $1$ \tabularnewline
\end{tabular}
\par\end{center}

\end{example}

\bigskip{}

\begin{example} \label{Ch7.Sec2.newE2} Multiplication group of the
AG-group in Example \ref{Ch7.Sec2.E2}. 

\begin{center}
\begin{tabular}{l|llllllll}
$\cdot$  & $L_{0}$  & $L_{1}$  & $L_{2}$  & $L_{3}$  & $R_{0}$  & $R_{1}$  & $R_{2}$  & $R_{3}$ \tabularnewline
\hline 
$L_{0}$  & $L_{0}$  & $L_{1}$  & $L_{2}$  & $L_{3}$  & $R_{0}$  & $R_{1}$  & $R_{2}$  & $R_{3}$ \tabularnewline
$L_{1}$  & $L_{1}$  & $L_{2}$  & $L_{3}$  & $L_{0}$  & $R_{3}$  & $R_{0}$  & $R_{1}$  & $R_{2}$ \tabularnewline
$L_{2}$  & $L_{2}$  & $L_{3}$  & $L_{0}$  & $L_{1}$  & $R_{2}$  & $R_{3}$  & $R_{0}$  & $R_{1}$ \tabularnewline
$L_{3}$  & $L_{3}$  & $L_{0}$  & $L_{1}$  & $L_{2}$  & $R_{1}$  & $R_{2}$  & $R_{3}$  & $R_{0}$ \tabularnewline
$R_{0}$  & $R_{0}$  & $R_{1}$  & $R_{2}$  & $R_{3}$  & $L_{0}$  & $L_{1}$  & $L_{2}$  & $L_{3}$ \tabularnewline
$R_{1}$  & $R_{1}$  & $R_{2}$  & $R_{3}$  & $R_{0}$  & $L_{3}$  & $L_{0}$  & $L_{1}$  & $L_{2}$ \tabularnewline
$R_{2}$  & $R_{2}$  & $R_{3}$  & $R_{0}$  & $R_{1}$  & $L_{1}$  & $L_{2}$  & $L_{3}$  & $L_{0}$ \tabularnewline
$R_{3}$  & $R_{3}$  & $R_{0}$  & $R_{1}$  & $R_{2}$  & $L_{2}$  & $L_{3}$  & $L_{0}$  & $L_{1}$ \tabularnewline
\end{tabular}
\par\end{center}

\end{example}

\bigskip{}

\noindent From Example \ref{Ch7.Sec2.newE2} we have the following
observations:
\begin{enumerate}
\item The multiplication group of an AG-group is not necessarily dihedral.
For example, $(L_{1}\cdot R_{3})^{2}=R_{2}^{2}=L_{3}\neq L_{0}$.
So here $M(G)$ is not $D_{4}$. 
\item From Examples \ref{Ch7.Sec2.newE1} and \ref{Ch7.Sec2.newE2} the
left sections in both the examples are $C_{3}$ and $C_{4}$ respectively. 
\end{enumerate}
\begin{theorem} Let $G$ be an AG-group. Let $a$ be an element of
$G$ distinct from $e.$ Then $a$ is self-inverse $\Longleftrightarrow R_{a}^{-1}=R_{a}$.
\end{theorem}

\begin{proof} Suppose $a$ is self-inverse. Since $R_{a}(x)=xa$,
then $R_{a}$ is of order $2$ , as $R_{a}(R_{a}(x))=(xa)a=(xa)a^{-1}=x\Longrightarrow$
$R_{a}^{2}=L_{e}\Longrightarrow R_{a}^{-1}=R_{a}.$

Conversely let $R_{a}^{2}=L_{e}$ then $R_{a}^{2}(x)=L_{e}(x)\,\forall x\in G\Longrightarrow(xa)a=ex=x.$
Now by left invertive law, $a^{2}x=x.$ This by right cancellation
implies $a^{2}=e$ or $a^{-1}=a.$ \end{proof}

\begin{remark} $R_{a}$ cannot fix all the elements of AG-group $G$.
For if we suppose that $R_{a}$ fixes all the elements. That is; $R_{a}(x)=x\,\forall x\in G\Longrightarrow xa=x\,\forall x\in G\Longrightarrow a$
is the right identity and hence $G$ is abelian. \end{remark}

\begin{theorem} For every AG-group $G$, the inner mapping group;
$Inn(G)=\left\{ L_{0},R_{0}\right\} $ is isomorphic to $C_{2}$.
\end{theorem}

\begin{proof} As $R_{a}(0)=0a=0.$ This implies that only $R_{0}$
maps $0$ on $0$. On the other hand $L_{0}(0)=0$ and no other $L_{a}$
can map $0$ on $0.$ Because let $L_{a}(0)=0$ where $a\neq0.$ Then
$a0=0.$ This implies $R_{0}(a)=0$. But $R_{0}(0)=0.$ This implies
that $R_{0}$ is not a permutation which is a contradiction. Hence
$Inn(G)=\left\{ L_{0},R_{0}\right\} \equiv C_{2}$. The following
table verifies the claim.

\[
\begin{tabular}{l|ll}
 \ensuremath{\cdot} &  \ensuremath{L_{0}} &  \ensuremath{R_{0}}\\
\hline \ensuremath{L_{0}} &  \ensuremath{L_{0}} &  \ensuremath{R_{0}}\\
\ensuremath{R_{0}} &  \ensuremath{R_{0}} &  \ensuremath{L_{0}}
\end{tabular}
\]

Hence the proof. \end{proof}

\noindent Again the following are some quick observations: \begin{enumerate}[(i)]

\item The $Inn(G)$ is not necessarily normal in $M(G)$ for example
consider the multiplication group of the AG-group given in \ref{Ch7.Sec2.E2}.
Here $L_{1}\left\{ L_{0},R_{0}\right\} =\left\{ L_{1},R_{3}\right\} \neq\left\{ L_{1},R_{1}\right\} =\left\{ L_{0},R_{0}\right\} L_{1}.$

\item For every AG-group $G$, $L_{S}$ being of index $2$ is normal
in $M(G)$ and hence $M(G)/L_{S}\equiv C_{2}.$ 

\item For every AG-group $G$, left multiplication group of $G$
coincides with $L_{S}$ and right multiplication group of $G$ coincides
with $M(G).$ 

\end{enumerate}

A non-associative quasigroup can be left distributive as well as right
distributive but a non-associative AG-group can neither be left distributive
nor right distributive as the following theorem shows.

\begin{theorem} \label{Ch7.Sec2.T8} Every left distributive AG-group
and every right distributive AG-group is abelian group. \end{theorem}

\begin{proof} Let $G$ be a left distributive AG-group. Then $\forall$
$a,b,c\in G,$ we have 
\begin{eqnarray*}
a(bc) & = & (ab)(ac)\\
 & = & (aa)(bc),\text{ by Lemma \ref{Pre-L1} (i)}\\
\Rightarrow\text{ }a & = & aa,\text{ by right cancellation}.
\end{eqnarray*}
This further implies that $G$ is an abelian group. The second part
is similar.\end{proof}

A non-associative quasigroup can be left distributive as well as right
distributive but a non-associative AG-group can neither be left distributive
nor right distributive as the following theorem shows.

\begin{theorem} \label{Ch7.Sec2.T13} If $G$ is an AG-group then
$M(G)$ cannot be the group of automorphisms of $L$. \end{theorem}

\begin{proof} Suppose on contrary that $M(G)$ is the group of automorphisms
of $G$. It means that every element of $M(G)$ is an automorphism
of $G$. Since $L_{a},R_{a}\in M(G)$ for all $a\in G$. Thus $L_{a}$
and $R_{a}$ are both automorphisms of $G$. So we can write
\begin{eqnarray*}
(xy)L_{a} & = & (x)L_{a}\cdot(y)L_{a},\textrm{ since }L_{a}\text{ is homomorphism}\\
\Rightarrow a(xy) & = & (ax)(ay)\text{ for all }x,y\in G
\end{eqnarray*}
Thus $G$ is left distributive. Similarly,
\begin{eqnarray*}
(xy)R_{a} & = & (x)R_{a}\cdot(y)R_{a},\textrm{ since }R_{a}\text{ is homomorphism}\\
\Rightarrow(xy)a & = & (xa)(ya)\text{ for all }x,y\in G
\end{eqnarray*}
Thus $G$ is right distributive. Hence $G$ is distributive, which
is a contradiction to Theorem \ref{Ch7.Sec2.T8}. Whence $M(G)$ of
an AG-group $G$ cannot be the group of automorphisms of $G$. \end{proof}

\begin{theorem} Let $e$ be the identity  and $x,y$ be any elements
of an AG-group $G$. Then:\end{theorem}

\begin{enumerate}[(i)]

\item $R_{x}^{-1}=R_{x^{-1}}$; 

\item $L_{x}^{-1}=L_{x^{-1}e}$. 

\end{enumerate}

\emph{Proof. }$(i)$ Since $G$ satisfies the right inverse property.
Therefore
\begin{eqnarray*}
(yx)x^{-1} & = & y\\
\Rightarrow R_{x^{-1}}R_{x}(y) & = & y=L_{e}(y)\,\forall x,y\in G\\
\Rightarrow R_{x^{-1}}R_{x} & = & L_{e}\Rightarrow R_{x}^{-1}=R_{x^{-1}}.
\end{eqnarray*}

$(ii)$ By Lemma \ref{Ch7.Sec2.L2} (iv)
\begin{eqnarray*}
L_{x}L_{x^{-1}e} & = & L_{(xe)(x^{-1}e)}=L_{(xx^{-1})e}=L_{e}\\
\Rightarrow L_{x}^{-1} & = & L_{x^{-1}e}.
\end{eqnarray*}
Hence the theorem.$\qed$\textbf{}

\end{document}